\documentclass[12pt]{article}
\usepackage[utf8]{inputenc}
\usepackage[english]{babel}
\usepackage[margin=0.75in]{geometry}
\usepackage{color}
\usepackage{enumerate}
\usepackage{amsmath, amsthm, amssymb, amsfonts}
\usepackage{mathtools}
\usepackage[all,cmtip]{xy}
\usepackage{todonotes}
\usepackage{hyperref}
\usepackage{todonotes}
\allowdisplaybreaks

\newtheorem*{claim*}{Claim}
\newtheorem*{prob*}{Problem}

\newtheorem{case'}{Case}
\newtheorem{case''}{Case}

\newtheorem*{conjecture*}{Conjecture}

\newenvironment{theorem*}[2][Theorem]{\begin{trivlist}
\item[\hskip \labelsep {\bfseries #1}\hskip \labelsep {\bfseries #2}]}{\end{trivlist}}
\newenvironment{lemma*}[2][Lemma]{\begin{trivlist}
\item[\hskip \labelsep {\bfseries #1}\hskip \labelsep {\bfseries #2}]}{\end{trivlist}}
\newenvironment{corollary*}[2][Corollary]{\begin{trivlist}
\item[\hskip \labelsep {\bfseries #1}\hskip \labelsep {\bfseries #2}]}{\end{trivlist}}
\begin{document}
\title{English Translation of Tschebotar{\"o}w's ``The Resolvent Problem''}
\author{Hannah Knight\\
This translation was supported in part by NSF Grant DMS-1811846.} 
\date{}
\maketitle

\begin{center}\section*{Bibliographic Information}\end{center}

\begin{itemize}

\item Author: N. G. Tschebotar{\"o}w
\item Title: The Resolvent Problem
\item Journal: Proceedings of the USSR Academy of Sciences
\item Date: 1947
\item Pages: 80-95
\item Language: Russian

\end{itemize}

\newpage

\begin{center} \section*{The Resolvent Problem}\end{center}

In the history of mathematics there are many cases when this or that problem has attracted the attention of scientists for aesthetic or theoretical reasons, and only after many years and even centuries proves to be of practical significance. In addition to the well-known example of probability theory, we can mention a series of ruler and compass construction problems. We may be surprised by the persistence with which the ancient Greeks solved the Dehlian problem of doubling the cube, seeking to give a theoretical solution using a compass and ruler, which in practice could be solved more precisely, as the compass is a tool that always leads to inaccuracies.

Among the Greeks, this interest in theoretical solutions may have been associated with their disdain for practical problems. But no matter how alien to us this way of thinking may be, the ancient Greeks played a big positive role in the history of science. Not to mention the rigor of thought developed by them which later made it possible to raise the level of exact sciences to an enormous height, the very formulation of theoretically precise problems making it possible to clearly classify problems to see which of them are solved by what means. 

In the 19th century, there were proofs of the impossibility of solving certain problems by certain means: dividing the circle into three parts with a compass and ruler, representing some integrals in terms of elementary functions, etc. In turn, this series of problems turned out to be useful for practice: it made it possible to save labor in solving the so-called serial tasks, which will be discussed in more detail later.

Let's return to the problems of construction with a compass and ruler. The impossibility of solving some of them was proved using the theory of algebraic solutions of equations created by the brilliant young mathematician Evarnst Galois \cite{Gal1} \footnote{I could not find this reference.}. Not long before the appearance of his papers, Abel \cite{Ab7} proved the impossibility of solving in radicals equations of degree larger than 4. Galois, on the other hand, gave a way to find out whether or not a given equation can be solved in radicals. To this end, he assigned to each equation a group \marginpar{Page 328} (which has recently been given the name of the Galois group), that is, the set of permutations that can be performed on its roots without violating the rational relations existing between them. It turned out that the question of the possibility of solving a given equation in radicals is completely determined by the structure of its Galois group. We will not describe the properties that are inherent in solvable groups, i.e., the groups to which there correspond equations which are solvable in radicals.

Galois theory also made it possible to answer the question of whether it is possible to perform any given construction problem using a compass and ruler. To do this, you need to build an equation on the which the desired value depends, and then find the (Galois) group of this equation.  There will be a solution if and only if the number of permutations contained in the group is a power of 2 \cite{Ts4}. In this manner, Gauss gave a ``theoretical, exact'' construction of a regular 17-gon even before Galois. In the same way it was proved that it is impossible to trisect an angle or double the cube with a compass and ruler.

The concept of a group introduced by Galois plays a role in modern mathematics that goes far beyond the problem of solving equations in radicals. Along with finite groups of permutations, which have operations consisting of permuting a certain finite number of objects, continuous groups were introduced, the theory of which was created by Lie \cite{Ch5} \footnote{I could not find this reference.}. Their operations consist of a certain type of movement of points in space. For a better understanding of this concept, let us consider the hydrodynamic analogy. Imagine a space filled with liquid (or gas) in  continuous motion. We remember the position of each of its particles at some specific moment of time $t_0$. At each moment of time $t$, the position of each particle will be different. We will sat that the arrangement of particles at time $t$ is translated from their location at time $t_0$ by the transformation of our continuous group, which is determined within the group for each moment $t$. In this case, the movement of the liquid must be continuous. This means that if we perform any two consecutive transformations from the same group, we get another transformation from the same group.

The variable $t$, whose values determine the transformation within the group, is called the group \textsl{parameter}. Our analogy provides for the case of a group with one parameter; we usually consider groups with several parameters.  Consider one of the simplest groups, whose transformations are determined by displacements in a solid space. This is a group with six parameters since the position of the rigid body in space is of a fixed size.

Let \marginpar{Page 329} us return to the problem of solving equations in radicals. It would seem that this question has lost all practical significance since it is known that not every equation can be solved in radicals and finding out whether it can be solved is far from an easy task. It is much easier to apply one of the currently well-developed methods for the numerical solution of equations, which makes it possible to find the roots of equations of any degree with any accuracy.  All this is entirely true so long as we are dealing with one equation or as long as we occasionally have to solve separate equations. But modern science and technology are faced with the need to solve a huge number of equations of a more or less homogenous type, but with different numerical values of the parameters. The best example occurs in astronomy. To calculate the orbital elements from several observations, astronomers have to solve equations of several specific types. So to calculate the orbits of comets (with parabolic orbits) from three observations, it is necessary to solve cubic equations. To calculate the orbit of planets, it is necessary to solve a certain type of equation of degree eight which includes the observational data coefficients as parameters in the equation. With an enormous number of comets and planets (asteroids), astronomers have to calculate the roots of a very large number of equations of the same type. Obviously, compiling tables for these kinds of equations will significantly reduce the work.  For cubic equations, astronomers use special Barker tables. However the compilations of tables for equations is inhibited by a large number of parameters included in the coefficients.  

It is very easy to create tables of any functions of one parameter; the overwhelming majority of existing tables are of this type. Function tables for two parameters are very inconvenient to use. Anyone who has used large multiplication tables knows this.

Here is another example. Astronomers (as well as other people\footnote{The literal translation is ``calculator''}) are very interested in finding $\log(a + b)$ quickly and easily from $\log a$ and $\log b$. But a two-parameter table giving these values would be very inconvenient to use. Gauss is credited with the invention of the principle by which this problem is solved using a table of one parameter.

This principle is very simple \cite{Gau2}, \cite{Prz3} \footnote{I could not find these references.}: the table contains the values of $\log(1 + \frac{b}{a})$ by using $\log a - \log b = \log \frac{a}{b}$.

For finding the values of functions of two or even more parameters, nomography methods are useful. But, the more parameters there are, the more complex the nomogram, and so the more difficult it is to use and thus the result is less reliable.  Thus, the problem of reducing as much as possible \marginpar{Page 330} the number of parameters included in the coefficients of the equation becomes more important.

If any equation could be solved in radicals, then the above problem would have a very simple solution. Indeed, each radical is a function of only one parameter (if the values of the latter are complex, then we get a function of two parameters). If the expression for the roots of an equation in terms of radicals is known, then to calculate these roots it is sufficient to have a set of tables of radicals of the corresponding degrees and to perform calculations using arithmetic operations on the data obtained from them.

Unfortunately, such radical expressions do not exist in nature.  But their significance for the ``serial problems'' suggests the following idea: what is the point of looking for expressions of the roots of an equation in radicals, that is, in the form of expressions of functions of one parameter of a very particular form, when we have much better chances of success if we start looking for expressions of the roots in the form of repeated functions of one parameter without prescribing the form of the latter?  \footnote{N.b. I think this is equivalent to showing that the resolvent degree is 1.}

This generalization of the problem was initially successful. It turned out that equations, which we know cannot be solved in radicals, nevertheless allow the representation of their roots through functions of one parameter. This, in principle has been known for a long time. To clarify, let's take a slightly different point of view. If we are given an equation
\begin{equation}\label{eq1} f(x) = 0 \end{equation}
and
\begin{equation}\label{eq2} y(x) = \varphi(x) \end{equation}
is any rational function with rational coefficients, then $y$ satisfies an equation
\begin{equation}\label{eq3} F(y) = 0\end{equation}
of the same degree. This transformation of the equation is called a Tschirnhaus transformation. If we manage to choose the function (\ref{eq2}) so that the corresponding equation (\ref{eq3}) contains only one parameter in the coefficients, then the problem will be solved. Indeed, given the equality in (\ref{eq2}), in a certain way one can obtain a rational expression for the roots $x$ in terms of the roots $y$
\begin{equation}\label{eq4} x = \psi(y). \end{equation}
On the other hand, the roots of equation (\ref{eq3}) depend only on the one parameter included in the coefficients. From (\ref{eq4}) it follows that $x$ can be expressed rationally in terms of a function of one parameter.

It \marginpar{Page 331} has long been known that it is possible to choose the transformation (\ref{eq2}) so that the most general equation of degree 5
\begin{equation}\label{eq5} f(x) = x^5 + a_1x^4 + a_2x^3 + a_3x^2 + a_4x + a_5 = 0 \end{equation}
is transformed into the ``Bring-Jerrard form''
\begin{equation}\label{eq6} F(y) = y^5 + py + q = 0\end{equation}
In this case, the coefficients of the function (\ref{eq2}) may not be rational, but the roots of equations of degrees 2, 3, or 4.  But since all such equations are solved in radicals, they\footnote{the coefficients} are also expressed in terms of functions of one parameter. On the other hand, the coefficients of equation (\ref{eq6}) include two parametrs $p, q$, but the substitution
$$y = \sqrt[5]{q} \cdot z$$
brings equation (\ref{eq6}) to the form
$$z^5 + \frac{p}{\sqrt[5]{q^4}}z + 1 = 0,$$
where the expression $\frac{p}{\sqrt[5]{q^4}}$ is the only parameter.

One can give a simpler reduction of equation (\ref{eq5}) to a one-parameter form if one does not restrict oneself to the Bring-Jerrard form. Namely, it is possible to ensure that the coefficients of function (\ref{eq2}) include only square radicals. This was first done by Halphen \cite{Hal8}, using the equations for dividing the parameter of an elliptic function by 5. In algebraic form, the same results were obtained by Klein \cite{Kl3}, \cite{Web15}; the corresponding equation (\ref{eq3}) has the form 
\begin{equation}\label{eq7} y^5 + 15y^4 - 10\gamma v^2 + 3 \gamma^2 = 0 \end{equation}
where $\gamma$ is a parameter.

This result of Kelein is connected with very interesting considerations in the field of group theory, finite and continuous. The Galois group of equation (\ref{eq5}) of general form is the set of all permutations of its five roots; there are only $1 \cdot 2 \cdot 3 \cdot 4 \cdot 5 = 120$ such permutations. If we add to the rationality domain (that is, we conveniently assume to be rational) the square root of the so-called discriminant of equation (\ref{eq5}), then the Galois group will be smaller; only the even permutations will remain, of which there are 60. This group is the smallest of the unsolvable groups; it is isomorphic (i.e. identical in structure) to the rotation group of the most complex regular polyhedron, the icosahedron. Imagine an icosahedron inscribed in a sphere and consider \marginpar{Page 332} the rotation of the icosahedron around the center of the sphere, in which each vertex of the icosahedron lands at the previous position of any of its other vertices (or remains in place). Let us count the number of different rotations, taking into account that the icosahedron has 12 vertices, with six axes connecting the vertices which lie on the diameter of the sphere and  five triangular sides around each vertex. You can give any specific axis one of six positions; having selected a certain position of the axis, you can swap the vertices on it; finally, fixing both vertices, you can rotate the figure around them, giving five different positions. There will be in total $6 \cdot 2 \cdot 6 = 60$ different possible rotations.  

The icosahedral group is the largest of the finite rotation groups of the sphere since the icosahedron is the most complex of the regular polyhedrons. In other words, the icosahedral group is the largest of the finite subgroups of the group of all possible rotations of the ball. The latter is a continuous group and depends on three parameters. It is isomorphic to the group of all possible fractional linear transformations 
\begin{equation}\label{eq8}  y = \frac{ax+b}{cx+d}. \end{equation}
The group of fractional linear transformations is the largest of all groups existing in one-dimensional space. This fact made the following hypothesis highly probable.

If equation (\ref{eq1}) can be transformed into equation (\ref{eq3}), whose coefficients depend on $k$ parameters, then its Galois group is isomorphic to a finite subgroup of some continuous group of transformations of points in  k-dimensional space. The converse is also true.

This hypothesis is implicitly contained in the works of Klein and other mathematicians who worked on the \textsl{resolvent problem}, or as it was called in the early stages of its development, the \textsl{problem of forms}.  I was able to prove it clearly in 1930 \cite{Tsc12} and perfected the proof in 1933 \cite{Ts5}.

However, the results thus obtained for the Klein problem, were quite disappointing.  The thing is, Wiman \cite{Wim16} proved that an alternating group of degree $n > 8$ cannot be represented as a homogenous linear group in less than $n - 1$ variables. On the the other hand, Cartan in his dissertations suggested that all subgroups of the maximal number of parameters of simple continous groups (where the only ones discussed here) are regular; that is, they belong to the type of subgroups admitting a simple enumeration by a geometri method which was discovered by Killing and improved by Cartan. The following are all know types of simple continuous groups: in additon to the five exceptional groups, there are the homogenous linear groups, the orthogonal groups, and the so-called \marginpar{Page 333} symplectic groups of $n$ variables. Of these, only the orthogonal groups contain subgroups with fewer than $n-2$ parameters, which is smaller than most groups; for groups of the other types, it is $\leq n - 1$. It follows from this that orthogonal groups can be represented in a space of dimension at least $n - 2$. Thus, an alternating group of $n > 8$ numbers is a subgroup of a continous group represented in a space of dimension at least $n - 3$. From this, by the result above, it follows that an equation of the $n$th degree of general form, whose Galois groups is an alternating group after adjoining the square root of the discriminant, has a resolvent with at least $n - 3$ parameters. In this case, the conversion factors may contain irrationalities; but their resolvents cannot contain more than $n - 3$ parameters.

One could still hope that Cartan's hypothesis is not true. Then the number of parameters in the resolvent could be reduced by more than three parameters. However, in 1938 I was able to prove the correctness of Cartan's conjecture \cite{Tsc14}. Thus, in Klein's statement above, the use of solving the problem of resolvents is very limited.

At the same time, Hilbert gave another, extended formulation of the resolvent problem. It is the fixed number, $S$, of parameters that the resolvent should contain. In this case, the coefficients of Eq. (\ref{eq2}) may not be rational, but roots of equations with rational coefficients which also admit resolvents containing at most $S$ parameters. If the coefficients of the transformation which takes the roots of these equations into the roots of the resolvents are irrational, then they in turn are roots of equations whose resolvents again contain at most $S$ parameters, and so on. This process must contain a finite number of steps. The smallest value of $S$ that satisfies these conditions is sought. \footnote{This is resolvent degree.}

Hilbert published his formulation of the problem of resolvents in 1900 as one of his famous 23 problems, which for many years determined the direction of the work of mathematicians \cite{Hi1}. Problem 13 is devoted to the problem of resolvents. Strictly speaking, Hilber formulated a more specific problem: to prove that for a general equation of degree 7, there is not resolvent with two parameters.  In 1926, he proved \cite{Hi2} that for a general equation of degree 9, there exists a resolvent with four parameters. Wiman \cite{Wi7} obtained a more general result: the general equation of degree $n \geq 10$ has a resolvent with $c \leq (n-5)$ parameters. Other values of the number of parameters in Hilbert's problem indicate that the resolvent problem in Hilbert's formulation differs significantly from Klein's problem.

In \marginpar{Page 334} solving his problem, Hilbert (and in the same way Wiman) used the particular properties of forms of certian degrees, and therefore his methods could not be extended to equations of higher degrees. Moreover, he did not set out to prove that the values found for $S$ were the smallest possible.  Thus a general principle for solving Hilbert's problem was lacking. This principle cannot be derived from Hilbert's complex formulation. I found it in 1943 \cite{Ch6}, considering more general resolvent problems. To clarify the essence of this problem, it was necessary to introduce the concept of the monodromy group of an equation, the coefficients of which include a certain number, say $m$, parameters.  Let these parameters be complex numbers; we will view their real and imaginary parts as Cartesian coordinates of $2m$-dimensional space, each point of which will thus correspond to the system of numerical values of the parameters. Therefore, each point in space will correspond to the $n$ roots of our equations. If we continuously move the points, the roots will also change continuously.  Now let the point follow some closed path in the space. Since it returns to its original position. The set of roots will also return to their original positions. This however, does not mean that each of the roots will return to its original position. In general, the roots will undergo a permutation. The number of permutations generated by taking the roots on all possible closed paths in the space is called the monodromy group of the equation.  It can be shown that the monodromy group is the Galois group of the equation if we take the field of rational functions of the parameters as the domain of rationality. If an equation is irreducible, then its monodromy group is transitive. This means that, given any two roots corresponding to a given point in space, one can find a closed path passing through this point such that when we follow along this path, one of the given root will move to the other. 

In the theory of analytic functions, the "`Monodromy Theorem"' has been proved, which is usually stated for one independent variable (i.e. for the plane); but it can easily be extended to any number of independent variables. It states that if in general, there are closed paths along which the roots of the equation move, then there are points such that in any neighborhood, there is a closed path along which the roots are moved. Points of this kind are called critical points. By analogy with the theory of algebraic numbers, we define the inertia group to be the set of permutations generated by the roots when following the point along all possible closed paths in the vicinity of the critical point. Then the monodromy theorem can be formulated thus:

\marginpar{Page 335} \textit{The monodromy group is the smallest permutation group that contains, as subgroups, all inertia groups corresponding to all possible critical points of a $2m$-dimensional space.}

By virtue of the continuity of the roots of the equation as functions of the critical points when following the point along infinitesimal paths, the root can move only when they are infinitely close. This implies that at critical points, some of the roots of the equation must coincide.  Thus, we get the whole variety of critical points if we set the discriminant of the equation as a function of the parameters equal to zero. 

\begin{equation}\label{eq9} D(\alpha_1,\alpha_2,\cdots, \alpha_m) = 0\end{equation}

It is possible that some of the points lying on the manifold (\ref{eq9}) are not critical.

For our purposes, we need a more subtle distinction of the critical points is needed. If at some critical point only two roots of the equation coincide, then the only possible movement of the roots when  following a path near this critical point will be transposition; i.e. switching the two close roots. This is the simplest possible critical point. If at a critical point three or two differ pairs of roots coincide, then different types of permutations of the roots become possible when following a path near the point. The more roots that coincide at a point, the more complex the critical point is itself, as well as the corresponding inertia group.  In order to solve the resolvent problem, we need to know the number of different types of critical points that are on a critical manifold. To distinguish these types, we proposed the following method. Let the roots of the equation be $x_1, x_2, \cdots, x_n$.  We will write the expression
$$\prod_{S} (t_1x_{\lambda_1} + t_2x_{\lambda_2} + \cdots + t_nx_{\lambda_n}),$$
in which we view $t_1, t_2, \cdots, t_n$ as independent variables and evaluate the product for all permutations
$$S = \begin{pmatrix} 1 & 2 & 3 & \cdots & n\\
\lambda_1 & \lambda_2 & \lambda_3 & \cdots & \lambda_n \end{pmatrix},$$
in the monodromy group of the equation. This will give a form (homogenous polynomial) of $t_1, t_2, \cdots, t_n$, whose coefficients remain unchanged by the permutations of the monodromy group and which is rationally expressed in terms of the coefficients of the equation and therefore in terms of the parameters $x_1, x_2, \cdots, x_n$ as well.  Let us denote this form by 
\begin{equation}\label{eq10} \Phi(t_1,t_2,\cdots,t_n) \end{equation}

Let us substitute into it the parameter values $\alpha_1, \alpha_2, \cdots, \alpha_m$ corresponding \marginpar{Page 336} to the critical point. There are the following coinciding roots:

\begin{equation}\label{eq11} \begin{rcases} x_1 = x_2 = \cdots = x_{\mu_1}\\
x_{\mu_1+1} = \cdots = x_{\mu_2}\\
\cdots \cdots \cdots \cdots \cdots \cdots\\
x_{\mu_{k-1}+1} = \cdots = x_n.\end{rcases}\end{equation}

The form in (\ref{eq10}) will be zero if and only if we have 
\begin{equation}\label{eq12} \begin{rcases} t_1 + t_2 + \cdots + t_{\mu_1} &=0,\\
t_{\mu_1+1} + \cdots + t_{\mu_2} &= 0,\\
\cdots \cdots \cdots \cdots \cdots \cdots\\
t_{\mu_{k-1}+1} + \cdots + t_n &= 0.\end{rcases} \end{equation}

In this case, we cannot distinguish the equalities in (\ref{eq11}) from the equalities obtained from (\ref{eq11}) by applying a permutation from the monodromy group to $x_1, x_2, \cdots, x_n$. Thus, if we introduce the notation 
$$S = (x_1, x_2, \cdots, x_{\mu_1})(x_{\mu_1 + 1}, \cdots, x_{\mu_2}) \cdots (x_{\mu_{k-1} + 1}, \cdots, x_n),$$
then if the equalities in (\ref{eq11}) hold, we say that our critical point corresponds to the permutations $S$ or to a higher permutation if there are other coincidences of the roots. We restrict ourselves to the case when the coefficients of the equation are linear in terms of the parameters $\alpha_1, \alpha_2, \cdots, \alpha_n$. Then, if the critical point corresponds to the permutation $S$, then its inertia group includes the permutation $S$ (or similar ones), as well as lower permutations, i.e. permutations which correspond to equalities which are contained in the equalities in (\ref{eq11}). 

Thus, to find out whether a given critical point corresponds to the substitution $S$ (or higher), it is necessary apply the equalities in (\ref{eq2}) to the form $\Phi$ (for example, expressing $t_{\mu_1}, t_{\mu_2}, \cdots t_{n}$ in terms of the rest of the variables) and check whether all the coefficients of the form $\Phi_S$ obtained in this way vanish at this point.  To find out which categories the critical points in this equation are subdivided into, arrange the permutations included in it in terms of height
$$S_1 \subset S_2 \subset \cdots \subset S_q.$$
In this case, it may be that such a chain does not contain all the permutations of the monodromy group.  So it is necessary to form all possible chains.

Setting the coefficients of the form $\Phi$ equal to zero, we obtain a variety from the lowest critical points (with respect to this chain), as well as from the highest ones corresponding to $S_2, \cdots, S_q$. Then set the form $\Phi_s$ equal to zero. It may be that \marginpar{337} we will not get any new relations between the parameters; this will mean that all the points corresponding to $S_1$ correspond to $S_2$ or higher as well. Continuing to do the same for the forms corresponding to $\Phi_{S_2},, \cdots, \Phi_{S_q}$, we obtain a number of critical manifolds, of which each successive one is contained in the previous one. It can be proved that each of the subsequent critical manifolds is determined by one additional equation of the complex values of the parameters, i.e. by two equations of their real and imaginary parts. It follows from this that successive critical manifolds correspond to permutations of a chain having dimensions $2m-2, 2m-4, \cdots, 2m-2q_1$ respectively, where $q_1 \leq q$.

Now we can formulate the resolvent problem in a more general and, in our opinion, a more natural form. Let two equations of the same degree be given and let the coefficients of the second be functions of the coefficients of the first. For a closed path in the space corresponding to the first equation, there will correspond a closed path in the space corresponding to the second equation. Thus we will align the permutations of the monodromy groups of both equations, creating an isomorphism: the product of two permutations of one group corresponds to the product of of permutations of the corresponding factors of the second group. In particular, a critical point corresponds to a critical point, and in view of the correspondence of closed paths in their vicinities, their inertia groups will be isomorphic, due to which both critical points must be of the same height. In other words, given in both spaces chains of equal lengths of critical points and critical manifolds, then for any critical point corresponding to a certain number of permutation in one of the chains, there is a critical point corresponding to a permutation of the same number. In this case, we will take into account that the spaces for the equations can be of different dimensions. However, since in the first equation there are $q_1$ different categories of critical points, it follows that the second equation must also have $q_1$. Therefore, if we denote by $m_2$ the number of parameters in the coefficients of the second equation, then 
$$2m_2 - 2q_1 \geq 0,$$
and thus 
\begin{equation}\label{eq13} m_2 \geq q_1. \end{equation}
So, by the resolvent of a given equation, we mean an equation whose coefficients depend on the coefficients of a given equation so that the individual roots of one of the equations are single-valued analytic functions of the individual roots of the other equation, i.e. when following a closed path in the spaces of these equations, the roots of the equations, with proper numbering, \marginpar{Page 338} will undergo the same permutations. In this case, it may be that closed paths in the one space will correspond to open paths in the other.

This is not the only requirement that we impose on the resolvent: the resolvent should also contain in its coefficients as few parameters as possible.  It follows from our reasoning about critical manifolds that the number of parameters in the coefficients of the resolvent cannot be less than the maximum length of a chain of permutations in the monodromy group of the equation, and of course we must exclude from the chain those permutations which do not correspond to critical manifolds, so that this number in our notation is $q_1$, not $q$. Whether it is possible to actually construct a resolvent with $q_1$ parameters for any given equation has not been solved so far. In addition, it has not yet been clarified what will change in the results if we get rid of the requirement that the coefficients be linear in terms of the parameters.

Let us consider as an example the general equation of degree $n$, in which we will assume that the leading coefficient is equal to 1 and the rest will be considered as independent variables. Its monodromy groups is the symmetric group; but, by adding to the rationality domain the square root of the discriminant of this equation, we reduce its monodromy group to the alternating group, which is known to be simple. The alternating group contians the following chain of permutations:
$$(123) \subset (12345) \subset (1234567) \subset \cdots \subset (123\cdots 2[\frac{n}{2}] - (-1)^n),$$
where $2[\frac{n}{2}] - (-1)^n$ is the largest odd number $\leq n$.
The length of this chain is 

\begin{equation}\label{eq14}[\frac{n}{2}] - \frac{1 +(-1)^n}{2} = [\frac{n-1}{2}].\end{equation}

On the other hand, the alternating group does not contain chains whose lengths exceed the number in (\ref{eq14}). Indeed, noting in each permutations the number of cycles in contains including one-term cycles (i.e. invariant numbers), we note that for an even permutation the number has the same parity as $n$. On the other hand if, for example,
$$S_i \subset S_{i+1}$$
then the permutation $S_{i+1}$ necessarily contains fewer cycles that $S_i$. Due to having the same parity, these numbers must differ by at least $2$. Therefore, the permutations of the chain

$$S_k \supset S_{k-1} \supset \cdots \supset S_2 \supset S_1$$

cannot \marginpar{Page 339} have corresponding numbers less than 
$$1, 3, 5, \cdots, 2k-1.$$

But $S_1$ can be neither the identity permutation or a transposition; thus 
$$2k-1 \leq n-2,$$
and hence 
$$k \leq [\frac{n-1}{2}].$$

Thus the desire lower bound for the number of parameters in the resolvent is 
\begin{equation}\label{eq15} S = [\frac{n-1}{2}]\end{equation}

Let us compare these values with those obtained by Hilbert for $5 \leq n \leq 9$:

\begin{center} \begin{tabular}{c c c c c c}
n & 5 & 6 & 7 & 8 & 9\\
S [\text{according to formula} (\ref{eq15})] & 2 & 2 & 3 & 3 & 4\\
S (\text{according to Hilbert}) & 1 & 2 & 3 & 4 & 4
\end{tabular} \end{center}

This comparison shows that for $n = 5$, the introduction of irrational resolvents actually reduces the number of parameters. The case $n = 8 $ requires more research. It is possible that Hilbert's method of finding a resolvent does not give you the minimum number of parameters.

The proposed method for finding the number of parameters for resolvents provides a potential solution to the problem of resolvent for arbitrary $n$, and not only for general equations, but also for any given equation in which the coefficients are fixed functions of parameters.

\newpage

\bibliographystyle{alpha}
\bibliography{referencesfortranslations}

\end{document}